%% file: kommutir_dissipativ03.tex
\documentclass[11pt,reqno]{amsart}
\usepackage{amssymb,amscd,amsbsy,mathrsfs}
\usepackage{mathdots}
\input{classstartscr}

\theoremstyle{remark}

\newtheorem*{rem*}{Remark}

\newcommand{\mE}{{\mathcal E}}

\newcommand{\OLA}{{\rm OL}_{\rm A}}

\newcommand{\COL}{{\rm CL}}

\newcommand{\CAo}{{\rm C}_{{\rm A},1}}
\newcommand{\Cb}{{\rm C}_{{\rm A},1}(\T^2)}
\newcommand{\CAbe}{{\rm C}_{{\rm A},\be}}
\newcommand{\CAb}{{\rm C}_{{\rm A},\be}(\C_+^2)}
\newcommand{\LIA}{\Li_{\rm A}}
\newcommand{\ri}{{\rm i}}

\newcommand\Li{{\rm Lip}}

\newcommand\dg{\frak D}

\begin{document}

\numberwithin{equation}{section}

\numberwithin{equation}{section}

\title{Functions of perturbed commuting dissipative operators}
\author{A.B. Aleksandrov and V.V. Peller}
\thanks{The research on \S\:4--\S\:6  is supported by 
Russian Science Foundation [grant number 18-11-00053].
The other results are partially supported by RFBR [grant number 20-01-00209a] 
and were conducted with the support of the
RUDN University Program 5-100}
\thanks{Corresponding author: V.V. Peller; email: peller@math.msu.edu}


\hfill To the memory of Hagen Neidhardt

\

\begin{abstract}
The main objective of the paper is to obtain sharp Lipschitz type estimates for the norm of operator differences $f(L_1,M_1)-f(L_2,M_2)$ for pairs $(L_1,M_1)$ and $(L_2,M_2)$ of commuting maximal dissipative operators. To obtain such estimates, we use double operator integrals with respect to semi-spectral measures associated with the pairs
$(L_1,M_1)$ and $(L_2,M_2)$. Note that the situation is considerably more complicated than in the case of functions of two commuting contractions and to overcome difficulties we had to elaborate new techniques. We deduce from the main result H\"older type estimates for operator differences as well as their estimates in Schatten--von Neumann norms.
\end{abstract} 

\maketitle

{\bf
\footnotesize
\tableofcontents
\normalsize
}

\setcounter{section}{0}
\section{\bf Introduction}
\setcounter{equation}{0}
\label{In}

\

We are going to study the behaviour of functions $f(L,M)$ of commuting maximal dissipative operators $L$ and $M$ under perturbation. In particular, we study conditions on functions $f$, under which the following Lipschitz type inequality holds:
$$
\|f(L_1,M_1)-f(L_2,M_2)\|\le\const\max\big\{\|L_2-L_1\|,\|M_2-M_1\|\big\}
$$
for commuting pairs $(L_1,L_2)$ and $(M_1,M_2)$ of maximal dissipative operators. We also obtain H\"older type estimates and estimates in Schatten--von Neumann norms 
$\bS_p$.

We extend earlier results obtained for functions of self-adjoint operators and unitary operators (see \cite{Pe1} and \cite{Pe3}), functions of contractions (see \cite{Pe3} and \cite{Pe4}), functions of dissipative operators (see \cite{AP3} and \cite{AP5}),
functions of normal operators (see \cite{APPS}) and functions of commuting contractions (see \cite{Pe6}).

Note that the methods developed in the papers mentioned above do not work in the case of commuting maximal dissipative operators and we had to elaborate new techniques.

In \S\:\ref{Disopfuis} we give a brief introduction in dissipative operators, we explain how to 
define the semi-spectral measure of such an operator and how to
construct a natural functional calculus. 

In \S\:\ref{KommuDis} we define a semi-spectral measure of a pair of commuting maximal dissipative operators $(L,M)$ and explain how to construct a natural functional calculus $f\mapsto f(L,M)$ for such a pair of operators. 

We give a brief introduction to Besov spaces in \S\:\ref{Besovy} and we define double operator integrals with respect to spectral and semi-spectral measure in \S\:\ref{DOI}.

In our paper \cite{AP5} we obtained Lipschitz type estimates for functions of maximal dissipative operators. We were unable to generalize the methods of \cite{AP5} to the case of functions of pairs of commuting dissipative operators. In \S\:\ref{odna} we offer a new approach in the case of functions of a single maximal dissipative operator. This new approach uses an idea, which will be used in \S\:\ref{dvepary} to obtain the main result of this paper. Note that such problems are closely related to the problem of finding an integrable spectral shift function for a pair of maximal dissipative operators with trace class difference. Recall that the last problem was solved in \cite{MNP1} and \cite{MNP2}, see also the earlier paper \cite{MN}.

In \S\:\ref{DdiHtp} we obtain estimates of divided differences in Haagerup tensor products. Such estimates will be used in \S\:\ref{dvepary}. Section \ref{nekommu} is devoted to functions of noncommuting maximal dissipative operators; this will also be used in \S\:\ref{dvepary}.

Finally, we deduce from the main result of the paper H\"older type estimates and estimates of operator differences in Schatten-von Neumann norms.

\

\section{\bf Dissipative operators and functional calculus}
\setcounter{equation}{0}
\label{Disopfuis}

\

In this section we give a brief introduction to maximal dissipative operators. We refer the reader to \cite{So}, \cite{AP3} and \cite{AP5} for more detailed information.

Recall that a not necessarily bounded linear operator $L$ in a Hilbert space $\h$ with domain $\cd_L$ is called {\it dissipative} if $\im(Lx,x)\ge0$ for every $x$ in $\cd_L$. It is called {\it maximal dissipative}  if it does not have a proper dissipative extension.

The {\it Cayley transform} of a dissipative operator $L$ is defined by
$$
T\df(L-{\rm i}I)(L+{\rm i}I)^{-1}
$$
with $\cd_T=(L+{\rm i}I)\cd_L$ and range $\Range T=(L-{\rm i}I)\cd_L$
($T$ is not densely defined in general). Then $T$ is a contraction, i.e., $\|Tu\|\le\|u\|$, $u\in \cd_T$, $1$ is not an eigenvalue of $T$, and $\Range(I-T)\df\{u-Tu:~u\in \cd_T\}$ is dense.

Conversely, if $T$ is a contraction defined on its domain $\cd_T$, $1$ is not an eigenvalue of $T$, and $\Range(I-T)$ is dense, then it is the Cayley transform of a dissipative operator $L$. Moreover, $L$ is the {\it inverse Cayley transform of} $T$:
$$
L={\rm i}(I+T)(I-T)^{-1},\quad \cd_L=\Range(I-T).
$$
Note that $L$ is maximal if and only if $\cd_T$ is the whole Hilbert space $\h$ .

Every dissipative operator has a maximal dissipative extension.
Every maximal dissipative operator $L$ is necessarily closed and
its spectrum $\s(L)$ is contained in the closed upper half-plane $\clos\C_+$.

If $L_1$ and $L_2$ are maximal dissipative operators, we say that the {\it operator $L_1-L_2$ is bounded} if there exists a bounded operator $K$ such that $L_2=L_1+K$.

Consider now the minimal unitary dilation $U$ of the Cayley transform $T$ of $L$, i.e., $U$ is a unitary operator defined on a Hilbert space $\K$ that contains $\h$ such that
$$
T^n=P_\h U^n\big|\h,\quad n\ge0,
$$
and $\K=\clos\spn\{U^nh:~h\in\h,~n\in\Z\}$. Since $1$ is not an eigenvalue of $T$, it follows that $1$ is not an eigenvalue of $U$ (see \cite{SNF}, Ch. II, \S\,6).
The Sz.-Nagy--Foia\c s functional calculus allows us to define a functional calculus for $T$ on the Banach algebra
$$
\CAo\df\big\{g\in H^\be:~g\quad
\mbox{is continuous on}\quad\T\setminus\{1\}~\big\}\!\!:
$$
$$
g(T)\df P_\h g(U)\Big|\h,\quad g\in\CAo.
$$
It is linear and multiplicative. Moreover,
$$
\|g(T)\|\le\|g\|_{H^\be},\quad g\in\CAo,
$$
(see \cite{SNF}, Ch. III).

We can define now a functional calculus for our dissipative operator on the
Banach algebra
$$
\CAbe=\big\{f\in H^\be(\C_+):~f\quad\mbox{is continuous on}\quad
\R\big\}:
$$
$$
f(L)\df \big(f\circ\o\big)(T),\quad f\in\CAbe,
$$
where $\o$ is the conformal map of $\dd$ onto $\C_+$ defined by
$\o(\z)\df{\rm i}(1+\z)(1-\z)^{-1}$, $\z\in\dd$.

We can extend now this functional calculus to the class $(z+\ri)\CAbe$. It is easy to see that this class contains the class $\LIA(\C_+)$ of of Lipschitz functions on $\clos\C_+$ that are analytic in $\C_+$.

Suppose that $f\in(z+\ri)\CAbe$.
We have
\bay
\label{fri}
f(\z)=\frac{f_\ri(\z)}{(\z+\ri)^{-1}},\quad\z\in\C_+,
\qquad
\mbox{where}\quad f_\ri(\z)\df\frac{f(\z)}{\z+\ri}.
\ey
Clearly, $f_\ri\in\CAbe$. The (possibly unbounded) operator $f(L)$ can be defined by
\bay
\label{f(L)}
f(L)\df(L+\ri I)f_\ri(L)
\ey
(see \cite{SNF}, Ch. IV, \S\;1). It follows from Theorem 1.1 of Ch. IV of \cite{SNF} that 
\bay
\label{fi(L)}
f(L)\supset f_\ri(L)(L+\ri I),
\ey
and so $\cd_{f(L)}\supset \cd_L$.

If $L$ is a maximal dissipative operator in a Hilbert space $\h$, we say that a self-adjoint operator $A$ in a Hilbert space
$\K$, $\K\supset\h$, is called a {\it resolvent self-adjoint dilation} of $L$
if
$$
(L-\l I)^{-1}=P_\h(A-\l I)^{-1}\Big|\h,\quad \im\l<0.
$$
The dilation is called {\it minimal} if
$$
\K=\clos\spn\big\{(A-\l I)^{-1}v:~v\in\h,~\im\l<0\big\}.
$$
If $f\in\CAbe$, then
$$
f(L)=P_\h f(A)\big|\h,\quad f\in\CAbe.
$$

A minimal resolvent self-adjoint dilation of a maximal dissipative operator always exists (and is unique up to a natural isomorphism). Indeed, it suffices to take a minimal unitary dilation of the Cayley transform of this operator and
apply the inverse Cayley transform to it. 

\medskip

{\bf Definition.} A {\it semi-spectral measure} $\mE$ on a $\s$-algebra $\Sigma$ of subsets of a set $\X$ is a map that takes values in the set of nonnegative operators on a Hilbert space $\h$ that is countably additive in the strong operator topology and such that $\mE(\X)=I$.

\medskip

By a theorem of Naimark (see \cite{Na}), a semi-spectral measure $\mE$ is a {\it compression} of a spectral measure, i.e., there is a Hilbert space $\K$ such that $\h\subset\K$, there is a spectral measure $E$ on $\K$ defined on the same $\s$-algebra $\Sigma$ such that 
$$
\mE(\D)=P_\h E(\D)\Big|\h\quad\D\in\Sigma.
$$
In this case $E$ is called a {\it dilation} of $\mE$. A dilation $E$ is called {\it minimal}
if 
$$
\K=\clos\spn\{E(\D)v:~\D\in\Sigma,~v\in\h\}.
$$
Note that it was shown in \cite{MM} that if $E$ is a minimal dilation of a semi-spectral measure $\mE$, then $\mE$ and $E$ are mutually absolutely continuous.

The semi-spectral measure $\mE_T$ of 
the Cayley transform $T$ of $L$  is defined by
$$
\mE_T(\D)\df P_\h E_U(\D)\Big|\h,
$$
where $\D$ is a Borel subset of $\T$ and $E_U$ is the spectral measure of the minimal unitary dilation $U$ of $T$.

Then
\bay
\label{ssu}
g(T)=\int_\T g(\z)\,d\mE_T(\z),\quad g\in\CAo.
\ey
The semi-spectral measure $\mE_L$ of $L$ is defined by
$$
\mE_L(\D)=\mE_T\big(\o^{-1}(\D)\big),\quad\D\quad\mbox{is a Borel subset of}
\quad\R.
$$
It follows easily from \rf{ssu} that
\bay
\label{fL}
f(L)=\int_\R f(x)\,d\mE_L(x),\quad f\in\CAbe.
\ey
If $A$ is the minimal self-adjoint resolvent dilation of $L$ and $E_A$ is the spectral measure of $A$, then
$$
\mE_L(\D)=P_\h E_A(\D)\big|\h,\quad\D\quad\mbox{is a Borel subset of}
\quad\R.
$$

\

\section{\bf Commuting dissipative operators and functional calculus}
\setcounter{equation}{0}
\label{KommuDis}

\

Recall that for not necessarily bounded linear operators $T$ and $R$ in a Hilbert space 
$\h$ with domains $\cd_T$ and $\cd_R$, the operator $TR$ is defined on its domain
$$
\cd_{TR}\df\{u\in\cd_R:\,\,\text{and}\,\, Ru\in\cd_T\}\quad\mbox{and}\quad TRu\df T(Ru),\quad
u\in\cd_{TR}.
$$

We are going to use repeatedly the following trivial observation.

\medskip

{\bf Remark.}
If $\M$ is a dense linear subset of $\h$ and $T$ is a closed linear operator with dense range such that $\M\subset\cd_T$, then $T(\M)$ is dense in $\h$. 
Indeed,
let $y\perp T(\M)$. Then for $x\in\M$, we have $0=(Tx,y)=(x,T^*y)$, and so $T^*y=0$,
i. e. $y\in(\Range T)^\perp$. 

\medskip

Let $L$ and $M$ be maximal dissipative operators in $\h$. We say that $L$ and $M$ {\it commute in the resolvent sense} if their resolvents $(L+\ri I)^{-1}$ and $(M+\ri I)^{-1}$ commute. Slightly abusing terminology, we say that  $L$ and $M$ commute if they commute in the resolvent sense. It is easy to verify that $(L+\ri I)^{-1}$ and $(M+\ri I)^{-1}$ commute if and only if the Cayley transforms $(L-\ri I)(L+\ri I)^{-1}$ and  $(M-\ri I)(M+\ri I)^{-1}$ of $L$ and $M$ commute.

\begin{lem}
\label{limi}
Let $L$ and $M$ be commuting maximal dissipative operators.
Put
$$
\cL\df\Range(L+\ri I)^{-1}(M+\ri I)^{-1}=\Range(M+\ri I)^{-1}(L+\ri I)^{-1}.
$$
Then $\cL$ is a dense subset of $\h$ and $\cL\subset\cd_L\cap\cd_M$.
\end{lem}

\Pf It follows trivially from the definition of $\cL$ that 
$\cL\subset\cd_{L+\ri I}=\cd_L$ and $\cL\subset\cd_{M+\ri I}=\cd_M$.
To show that $\cL$ is dense in $\h$ it suffices to observe that 
$$
\cL= (L+\ri I)^{-1}(\Range(M+\ri I)^{-1})=(L+\ri I)^{-1}\cd_M,
$$
$\cd_M$ is dense in $\h$ and $\Range(L+\ri I)^{-1}$ is dense in $\h$. $\bl$

\medskip

%
%
%

Note that it follows from the definition of $\cL$ that
$$
(L+\ri I)\cL\subset \cd_M=\cd_{M+\ri I}\quad\mbox{and}\quad
(M+\ri I)\cL\subset\cd_L=\cd_{L+\ri I}.
$$

Lemma \ref{limi} implies that 
\bay
\label{lm}
\cd_{(L+\ri I)(M+\ri I)}=\cd_{(M+\ri I)(L+\ri I)}=\cL
\quad\mbox{and}\quad
(L+\ri I)(M+\ri I)=(M+\ri I)(L+\ri I).
\ey.

It is easy to see that $\cd_{LM}\ne\cd_{ML}$ in general.
Indeed, suppose that $L$ is an unbounded maximal dissipative operator with bounded $L^{-1}$.
Put $M=-L^{-1}$. Then $\cd_L\ne\h$ and $\cd_M=\h$. Hence, $\cd_{LM}=\cd_M=\h$
but $\cd_{ML}=\cd_L\ne\h$.


\begin{lem}
\label{LMML}
Let $L$ and $M$ be commuting maximal dissipative operators and let
$u\in\cd_L\cap\cd_M$. The following are equivalent:

{\rm(a)} $u\in\cd_{LM}$;

{\rm(b)} $u\in\cd_{ML}$;

{\rm(c)} $u\in\cL$.

\noindent
Moreover, $L(Mu)=M(Lu)$ if  one of statements {\rm(a)}, {\rm(b)} or {\rm(c)} holds.
\end{lem}

\Pf The result follows from \rf{lm}. Indeed, let $u\in\cd_L\cap\cd_M$, then $u\in\cd_{(L+\ri I)(M+\ri I)}$ if and only if
$u\in\cd_{LM}$ and $u\in\cd_{(M+\ri I)(L+\ri I)}$
if and only if
$u\in\cd_{ML}$. Moreover,  $u\in\cL$ if and only $u\in\cd_{(L+\ri I)(M+\ri I)}$.
It remains to observe that $u\in\cd_{(L+\ri I)(M+\ri I)}$ if and only if $u\in\cd_{LM}$. Finally, the equality $LMu=MLu$ trivially follows from
the equality $(L+\ri I)(M+\ri I)u=(M+\ri I)(L+\ri I)u$. $\bl$

\begin{lem}
\label{LM33}
Let $L$ and $M$ be maximal dissipative operators.
Suppose that $\M$ is a dense linear subset of $\h$ and
$\M\subset\cd_{LM}\cap\cd_{ML}$. If
$L(Mu)=M(Lu)$ for all $u\in\M$, then $L$ and $M$ commute.
\end{lem}

\Pf Let us first show that $(L+\ri I)(M+\ri I)\M$ is dense in $\h$.  Clearly,
\bay
\label{comLM}
(M+\ri I)^{-1}(L+\ri I)^{-1}u=(L+\ri I)^{-1}(M+\ri I)^{-1}u
\ey
for all $u\in(L+\ri I)(M+\ri I)\M=(M+\ri I)(L+\ri I)\M$.
Hence, \rf{comLM} holds for all $u\in\h$.
$\bl$

\begin{thm} Let $L$ and $M$ be maximal dissipative operators.
Then the following are equivalent:

{\rm(a)} $L$ and $M$ commute (in the resolvent sense);

{\rm(b)}  $\cd_{LM}\cap\cd_{ML}$ is dense in $\h$ and
$L(Mu)=M(Lu)$ for all $u\in\cd_{LM}\cap\cd_{ML}$;

{\rm(c)} $(L+\ri I)(M+\ri I)=(M+\ri I)(L+\ri I)$;

{\rm(d)} the Cayley transforms of $L$ and $M$ commute.
\end{thm}

\Pf The implication $(\rm a)\Rightarrow(\rm b)$ follows from Lemma \ref{LMML}.
The implication $(\rm b)\Rightarrow(\rm a)$ follows from Lemma \ref{LM33}. 
Next, (c) follows from (a) by \rf{lm}. Finally, let us establish that 
$(\rm c)\Rightarrow(\rm a)$. Suppose that
$$
(L+\ri I)(M+\ri L)=(M+\ri L)(L+\ri I).
$$
Then
$$
(L+\ri I)(M+\ri L)(L+\ri I)^{-1}(M+\ri I)^{-1}=I,
$$
whence
$$
(L+\ri I)^{-1}(M+\ri I)^{-1}=(M+\ri I)^{-1}(L+\ri I)^{-1}.
$$
The equivalence of (a) and (d) was mentioned above. $\bl$

\medskip

Consider now the Cayley transforms
$T=(L-\ri I)(L+\ri I)^{-1}$ and  $R=(M-\ri I)(M+\ri I)^{-1}$
of $L$ and $M$.
It is well known that $\Ker(I-T)=\Ker(I-T^*)=\Ker(I-R)=\Ker(I-R^*)=\{0\}$.

By Ando's theorem \cite{An}, $T$ and $R$ have commuting unitary dilations. In other words, there exist a Hilbert space $\K$ that contains $\h$ and commuting unitary operators $U$ and $V$ on $\K$ such that
$$
P_\h U^jV^k\big|\h=T^jR^k,\quad j,~k\ge0,
$$
where $P_\h$ is the orthogonal projections onto $\h$. 

Let us show that without loss of generality we may assume that 
$1\not\in\s_{\rm p}(U)$ and $1\not\in\s_{\rm p}(V)$.
Indeed, let $\N\df\Ker(U-I)$. Then $\N$ is a reducing subspace of $U$. Let us show that $\N$ also reduces $V$. Indeed, suppose that $x\in\N$.
Then
$$
UVx=VUx=Vx\quad\mbox{and}\quad UV^{-1}x=V^{-1}Ux=V^{-1}x
$$
and so both $Vx$ and $V^{-1}x$ belong to $\N$.

It is well known that the minimal unitary dilation of $T$ has no point spectrum at 1, and so $\N\perp\h$. We can replace now $\K$ with $\K\ominus\N$ and the unitary operators $U$ and $V$ with their restrictions to $\K\ominus\N$. Thus, we may assume that 1 is not an eigenvalue of $U$. Similarly, we can consider $\Ker(V-I)$ and consider the restrictions of the unitary dilations to its orthogonal complement.

We proceed now to the construction of a functional calculus
$$
f\mapsto f(L,M).
$$
First, we construct a functional calculus for the Cayley transforms $T$ and $R$.

We define the subclass $\CAo(\T^2)$ of the Hardy class $H^\be(\T^2)$ of functions on the bidisk $\T^2$ by
$$
\Cb\df\big\{g\in H^\be(\T^2):~g\quad
\mbox{is continuous on}\quad\T^2\setminus\{1,1\}~\big\}.
$$
For a function $g$ in $\Cb$, we put
$$
g(T,R)\df P_\h g(U,V)\Big|\h.
$$
Note that $g(T,R)$ does not depend on the choice of commuting unitary dilations. Indeed, it is easy to see that
$$
g(T,R)x=\lim_{r\to1}g_r(T,R)x,\quad x\in\h,
$$
where $g_r(\z,\t)\df g(r\z,r\t)$,
and
$$
g_r(T,R)=\sum_{j,k\ge0}r^{j+k}\widehat g(j,k)T^jR^k
$$
which does not depend on $U$ or $V$.

Consider now the class 
$$
\CAb\df\big\{f\in H^\be(\C_+^2):~f\quad\mbox{is continuous on}\quad
\R^2\big\},
$$
where $H^\be(\C_+^2)$ is the Hardy class of bounded analytic functions on $\C_+^2$.
We can define now the functional calculus $f\mapsto f(L,M)$ on the class 
$\CAb$ by
$$
f(L,M)\df(f\circ\bs{\o})(T,R),
$$
where 
$$
\bs{\o}(\z_1,\z_2)\df
\big({\rm i}(1+\z_1)(1-\z_1)^{-1},{\rm i}(1+\z_2)(1-\z_2)^{-1}\big),
\quad\z_1,~\z_2\in\dd.
$$

Consider now the joint spectral measure $E_{U,V}$ of the commuting pair $(U,V)$. 
We can define now the semi-spectral measure $\mE$ 
by
$$
\mE(\D)\df P_\h E_{U,V}(\D)\big|\h,
$$ 
where $\D$ is a Borel subset of $\T^2$. We say that $\mE$ is {\it a semi-spectral measure of the pair $(T,R)$ of commuting contractions}. We do not know know whether it can depend on the choice of the commuting unitary dilation $(U,V)$. However, we can select one of such semi-spectral measures and denote it by $\mE_{T,R}$. It is easy to see that
$$
g(T,R)=\int_{\T^2}g(\z_1,\z_2)\,d\mE_{T,R}(\z_1,\z_2),\quad g\in\Cb.
$$

This allows us to define the semi-spectral measure
$\mE_{L,M}$ of the pair $(L,M)$  by
$$
\mE_{L,M}(\L)=\mE_{U,V}(\bs{\o}^{-1}(\L)),
$$
where $\L$ is a Borel subset of $\R^2$. It is easy to verify that 
$$
f(L,M)=\int_{\R^2}f(\z)\,d\mE_{L,M}(\z)
$$
for an arbitrary function $f$ in $\CAb$. We are going to say that $\mE_{L,M}$ is {\it a semi-spectral measure of the commuting pair $(L,M)$}.

We can extend now this functional calculus to the class $(z_1+\ri)(z_2+\ri)\CAb$,
which contains  the class $\LIA(\C_+^2)$  of Lipschitz
functions on $(\clos\C_+)^2$ that are analytic in $\C_+^2$. 

If $f$ belongs to this class, we put
$$
f_\ri(\z_1,\z_2)\df\frac{f(\z_1,\z_2)}{(\z_1+\ri)(\z_2+\ri)},\quad\z_1,~\z_2\in\C_+.
$$
Clearly, $f_\ri\in\CAb$.  We define the (not necessarily bounded operator) $f(L,M)$ by
\bay
\label{f(L,M)}
f(L,M)=(L+\ri I)(M+\ri I)f_\ri(L,M).
\ey
As in the case of functions of one dissipative operator, we have
$$
f(L,M)\supset f_\ri(L,M)(L+\ri I)(M+\ri I)
$$
(cf. \rf{fi(L)}).

To see this, we can use the same argument as in the case of functions of a single maximal dissipative operator (see Ch. IV of \cite{SNF}). We present a proof for the sake of reader's convenience.

Equality \rf{f(L,M)} implies
\bey
f(L,M)((L+\ri I)(M+\ri I))^{-1}=(L+\ri I)(M+\ri I)f_\ri(L,M)((L+\ri I)(M+\ri I))^{-1}\\
=(L+\ri I)(M+\ri I)((L+\ri I)(M+\ri I))^{-1}f_\ri(L,M)=f_\ri(L,M).
\eey
Hence,
$$
f(L,M)((L+\ri I)(M+\ri I))^{-1}(L+\ri I)(M+\ri I)=f_\ri(L,M)(L+\ri I)(M+\ri I)
$$
Taking into account that $((L+\ri I)(M+\ri I))^{-1}(L+\ri I)(M+\ri I)\subset I$,
we see that
$$
f(L,M)\supset f_\ri(L,M)(L+\ri I)(M+\ri I).
$$

\

\section{\bf Besov spaces}
\setcounter{equation}{0}
\label{Besovy}

\

In this paper we deal with the homogeneous Besov class $B_{\be,1}^1(\R^2)$. Moreover, we need only the
analytic Besov class $\big(B_{\be,1}^1\big)_+(\R^2)$.
We refer the reader to the book \cite{Pee} and the papers \cite{ANP} and \cite{AP4} for more information on
Besov classes $B_{p,q}^s(\R^d)$. Here we give the definition only in the case
when $p=\be$, $q=s=1$.
Let $w$ be an infinitely differentiable function on $\R$ such
that
\bay
\label{w}
w\ge0,\quad\supp w\subset\left[\frac12,2\right],\quad\mbox{and} \quad w(t)=1-w\left(\frac t2\right)\quad\mbox{for}\quad t\in[1,2].
\ey

Consider the functions $W_n$, $n\in\Z$, on $\R^d$ such that 
$$
\big(\F W_n\big)(x)=w\left(\frac{\|x\|_2}{2^n}\right),\quad n\in\Z, \quad x=(x_1,\cdots,x_d),
\quad\|x\|_2\df\left(\sum_{j=1}^dx_j^2\right)^{1/2},
$$
where $\F$ is the {\it Fourier transform} defined on $L^1\big(\R^d\big)$ by
$$
\big(\F f\big)(t)=\!\int\limits_{\R^d} f(x)e^{-{\rm i}(x,t)}\,dx,\!\quad 
x=(x_1,\cdots,x_d),
\quad t=(t_1,\cdots,t_d), \!\quad(x,t)\df \sum_{j=1}^dx_jt_j.
$$
Clearly,
$$
\sum_{n\in\Z}(\F W_n)(t)=1,\quad t\in\R^d\setminus\{0\}.
$$

With each tempered distribution $f\in{\mathscr S}^\prime\big(\R^d\big)$, we
associate the sequence $\{f_n\}_{n\in\Z}$,
\bay
\label{fn}
f_n\df f*W_n.
\ey
The formal series
$
\sum_{n\in\Z}f_n
$
is a Littlewood--Paley type expansion of $f$. 
This series does not necessarily converge to $f$. 

Initially we define the (homogeneous) Besov class $\dot B^1_{\be,1}\big(\R^d\big)$ as the space of 
$f\in{\mathscr S}^\prime(\R^n)$
such that
\bay
\label{<be}
\|f\|_{B^1_{\be,1}}\df\sum_{n\in\Z}2^n\|f_n\|_{L^\be}<\be.
\ey
According to this definition, the space $\dot B^1_{\be,1}(\R^n)$ contains all polynomials
and all polynomials $f$ satisfy the equality $\|f\|_{B^1_{\be,1}}=0$. Moreover, the distribution $f$ is determined by the sequence $\{f_n\}_{n\in\Z}$
uniquely up to a polynomial. It is easy to see that the series 
$\sum_{n\ge0}f_n$ converges in ${\mathscr S}^\prime(\R^d)$.
However, the series $\sum_{n<0}f_n$ can diverge in general. Obviously, the series
\bay
\label{ryad}
\sum_{n<0}\frac{\partial f_n}{\partial x_j},\quad
1\le j\le d,
\ey
converges uniformly on $\R^d$. 
Now we say that  $f$ belongs to the {\it homogeneous Besov class} $B^1_{\be,1}(\R^d)$ if \rf{<be} holds and 
$$
\frac{\partial f}{\partial x_j}=\sum_{n\in\Z}\frac{\partial f_n}{\partial x_j},\quad
1\le j\le d.
$$

A function $f$ is determined uniquely by the sequence $\{f_n\}_{n\in\Z}$ up
to a a constant, and a polynomial $g$ belongs to 
$B^1_{\be,1}\big(\R^d\big)$
if and only if it is a constant.

Put
$$
\big(B_{\be,1}^1\big)_+(\R^d)\df\big\{f\in B_{\be,1}^1(\R^d):\supp\F f\subset[0,+\be)^d\big\}.
$$
Functions of class $B^1_{\be,1}\big(\R^d\big)$ can be considered as functions 
analytic in $\C_+^d$ and continuous on $\clos\C_+^d$.

\

\section{\bf Double operator integrals}
\setcounter{equation}{0}
\label{DOI}

\

Double operator integrals with respect to spectral measures are expressions of the form
\bay
\label{doiQE1E2}
\iint_{\X_1\times\X_2}\Phi(x_1,x_2)\,dE_1(x)Q\,dE_2(y).
\ey
Here $E_1$ and $E_2$ are spectral measures on a Hilbert space, 
$\Phi$ is a bounded measurable function and $Q$ is a bounded linear operator on $\h$. 

Double operator integrals appeared first in \cite{DK}. Later Birman and Solomyak elaborated a beautiful theory of double operator integrals in \cite{BS1}--\cite{BS3}.
The starting point of the Birman--Solomyak approach is the definition of the double operator integral in the case when $Q$ belongs to the Hilbert--Schmidt class $\bS_2$,
in which case the double operator integral in \rf{doiQE1E2}
must also be in $\bS_2$ and 
$\left\|\iint\Phi\,dE_1Q\,dE_2\right\|_{\bS_2}\le\sup|\Phi|\cdot\|Q\|_{\bS_2}$.

However, if we want to define the double operator integral  for an arbitrary bounded linear operator $Q$, we have to impose additional assumptions on $\Phi$. The natural class of functions $\Phi$, for which the double operator integral determines a bounded linear operator for all bounded $Q$ is the {\it Haagerup tensor product} 
$L^\be_{E_1}\otimes_{\rm h}L^\be_{E_2}$, which consists of functions $\Phi$ of the form
$$
\Phi(x_1,x_2)=\sum_n\f_n(x_1)\psi_n(x_2),
$$
where $\f_n$ and $\psi_n$ are functions in $L^\be_{E_1}$ and $L^\be_{E_2}$ satisfying
$$
\left\|\sum_n|\f_n|^2\right\|_{L^\be_{E_1}}\left\|\sum_n|\psi_n|^2\right\|_{L^\be_{E_2}}<\be
$$
(see \cite{BS3}, \cite{Pe1} and \cite {AP4}). In this case
$$
\iint\Phi\,dE_1Q\,dE_2=\sum_n\left(\int\f_n\,dE_1\right)Q\left(\int\psi_n\,dE_2\right)
$$
and the series converges in the weak operator topology.

In a similar way we can define Haagerup tensor products of subspaces of $L^\be$.

Double operator integrals with respect to semi-spectral measures were first considered in
\cite{Pe2}.  Let $\mE_1$ and $\mE_2$ be semi-spectral measures on a Hilbert space $\h$. 
Let $E_1$ and $E_2$ are spectral measures on Hilbert spaces $\K_1$ and $\K_2$, which are minimal dilations of $\mE_1$ and $\mE_2$ (recall that $\K_1\supset\h$ and 
$\K_2\supset\h$). As we have mentioned in \S\:\ref{Disopfuis}, $E_j$ and $\mE_j$
are mutually absolutely continuous, $j=1,2$. 

Suppose that $\Phi$ is a function in the Haagerup tensor product 
$L^\be_{E_1}\otimes_{\rm h}L^\be_{E_2}=L^\be_{\mE_1}\otimes_{\rm h}L^\be_{\mE_2}$
and let $Q$ be a bounded linear operator on $\h$. Then
$$
\iint\Phi(x_1,x_2)\,d\mE_1(x_1)Q\,d\mE_2(x_2)\df 
P^{\K_1}_\h\iint\Phi(x_1,x_2)\,dE_1(x_1)Q_\flat\,dE_2(x_2)\Big|\h,
$$
where $Q_\flat$ is the bounded linear from $\K_2$ to $\K_1$ defined by 
$Q_\flat u=QP^{\K_2}_\h$, $u\in\K_2$, where $P^{\K_j}_\h$ 
is the orthogonal projections from $\K_j$ onto $\h$, $j=1,2$.

\

\section{\bf A new approach to operator Lipschitz estimates\\ for functions of a single dissipative operator}
\setcounter{equation}{0}
\label{odna}

\

The purpose of this section is to give an alternative proof of Theorem \ref{OLB9} below, which is the main result of \cite{AP4}. Recall that a function $f$ continuous on the  upper half-plane $\clos\C_+$ is called an {\it operator Lipschitz function on $\clos\C_+$}
 if there exists a constant $C$ such that
\bay
\label{lipN}
\|f(N_1)-f(N_2)\|\le C\|N_1-N_2\|
\ey
for arbitrary normal operators $N_1$ and $N_2$ with spectra in $\clos\C_+$.
We denote by $\OLA(\C_+)$ the space of operator Lipschitz functions on $\clos\C_+$ that
are analytic in $\C_+$. For $f\in\OLA(\C_+)$, we denote by $\|f\|_{\OLA(\C_+)}$ 
be the smallest constant C, for which \rf{lipN} holds. Note that the class $\OLA(\C_+)$ coincides with the class $\COL(\C_+)$ of commutator Lipschitz functions on $\clos\C_+$, see \cite{AP4}. Let us also mention that the class $\OLA(\C_+)$ can be identified with the class operator of Lipschitz functions on $\R$ whose derivative belongs to $H^\be(\C_+)$.

Suppose that $f\in\OLA(\C_+)$ and consider its restriction to the real line $\R$. We are going to use the same notation $f$ for the restriction. Then $f$ is an {\it operator Lipschitz function on the real line} $\R$, i.e.,
$$
\|f(A)-f(B)\|\le\const\|A-B\|
$$
whenever $A$ and $B$ are self-adjoint operators with bounded $A-B$. By a theorem of Johnson and Williams (see \cite{AP4}), such a function $f$ must be differentiable everywhere on $\R$. Thus, we can define the divided difference $\dg f$ 
$$
(\dg f)(x,y)\df\left\{\begin{array}{ll}\frac{f(x)-f(y)}{x-y},&x\ne y,\\[.2cm]
0,&x=y.
\end{array}\right.
$$
It is well known (see \cite{AP4}, Th. 3.9.6) that for operator Lipschitz functions $f$ on $\R$, the divided difference $\dg f$ is a Schur multiplier with respect to arbitrary Borel spectral measures $E_1$ and $E_2$, and so
$$
\left\|\iint_{\R\times\R}(\dg f)(x,y)\,dE_1(x)Q\,dE_2(y)\right\|\le\const\|Q\|
$$
for an arbitrary bounded operator $Q$. Clearly, the same estimate holds if we replace spectral measures with semi-spectral measures.

\begin{thm}
\label{OLB9}
Let $f\in\OLA(\C_+)$. Then for arbitrary maximal dissipative operators
$L$ and $M$ with bounded $L-M$, the following formula holds:
\bay
\label{BSd9}
f(L)-f(M)=\iint_{\R\times\R}(\dg f)(x,y)\,d\mE_L(x)(L-M)\,d\mE_M(y),
\ey
and so
$$
\|f(L)-f(M)\|\le\const\|f\|_{\OLA(\C_+)}\|L-M\|.
$$
\end{thm}

{\bf Comment.} We assume that the operators $L-M$ and $f(L)-f(M)$ are extended by continuity to the whole space. 

\medskip

To prove the theorem, we need a lemma.
Put $\iota(z)\df z(1-\ri z)^{-1}$,
$L_1\df\iota(L)$ and $M_1\df\iota(M)$.

\begin{lem}
\label{LeMe9}
Let $L$,  $M$ and $f$ satisfy the hypotheses of Theorem {\em\ref{OLB9}}.
Then
\begin{multline}
\label{L1M1Df}
\iint_{\R\times\R}(\dg f)(x,y)\,d\mE_L(x)(L_1-M_1)\,d\mE_M(y)\\
=(I-\ri L)^{-1}\left(\iint_{\R\times\R}(\dg f)(x,y)\,d\mE_L(x)(L-M)\,d\mE_M(y)\right)(I-\ri M)^{-1}.
\end{multline}
\end{lem}

\Pf We have
\bey
\label{Hilb9}
L_1-M_1=
L(I-\ri L)^{-1}-M(I-\ri M)^{-1}
=(I-\ri L)^{-1}(L-M)(I-\ri M)^{-1}
\eey
by a version of the Hilbert resolvent identity.
Applying this identity, we obtain the desired equality. $\bl$

\medskip

{\bf Proof of Theorem \ref{OLB9}.} 
In view of \rf{L1M1Df}, it suffices to prove that
\bay
\label{pred9}
\iint\limits_{\R\times\R}(\dg f)(x,y)\,d\mE_L(x)(L_1-M_1)\,d\mE_M(y)=(I-\ri L)^{-1}(f(L)-f(M))(I-\ri M)^{-1}.
\ey
Let us first assume that $f$ is bounded. We have
\begin{align*}
\iint\limits_{\R\times\R}&(\dg f)(x,y)\,d\mE_L(x)(L_1-M_1)\,d\mE_M(y)=
\iint\limits_{\R\times\R}(\dg f)(x,y)\,d\mE_L(x)L_1\,d\mE_M(y)\\[.2cm]
&-\iint\limits_{\R\times\R}(\dg f)(x,y)\,d\mE_L(x)M_1\,d\mE_M(y)=
\iint\limits_{\R\times\R}\iota(x)(\dg f)(x,y)\,d\mE_L(x)\,d\mE_M(y)\\[.2cm]
&-\iint\limits_{\R\times\R}\iota(y)(\dg f)(x,y)\,d\mE_L(x)\,d\mE_M(y)=
\iint\limits_{\R\times\R}(f(x)-f(y))(\dg \iota)(x,y)\,d\mE_L(x)\,d\mE_M(y)\\[.2cm]
&=(I-\ri L)^{-1}\left(\,\,\iint\limits_{\R\times\R}(f(x)-f(y))\,d\mE_L(x)\,d\mE_M(y)\right)(I-\ri M)^{-1}\\[.2cm]
&=(I-\ri L)^{-1}(f(L)-f(M))(I-\ri M)^{-1}
\end{align*}
which proves \rf{pred9}.

Let now $f$ be an arbitrary function in $\OLA(\C_+)$. Put $f_\e(z)\df(1-\ri\e z)^{-1}f(z)$, where $\e>0$.
Clearly, the functions $f_\e$  are bounded. Moreover, $f_\e(z)\in\OLA(\C_+)$.
This follows from Theorem 4.6 in \cite{A}

We prove this for the reader's convenience. 
It suffices to consider the case when $f(0)=0$.
Then $f$ can be represented in the form $f(z)=zg(z)$, where $g\in H^\be\cap{\rm C}(\R)$. We have
\begin{align*}
f_\e(N_1)-f_\e(N_2)&=(I-\ri\e N_1)^{-1}(f(N_1)-f(N_2))\\[.2cm]
&+((I-\ri\e N_1)^{-1}-(I-\ri\e N_2)^{-1})f(N_2)\\[.2cm]
&=(I-\ri\e N_1)^{-1}(f(N_1)-f(N_2))-(I-\ri\e N_1)^{-1}(N_1-N_2)g(N_2)\\[.2cm]
&+(N_1(I-\ri\e N_1)^{-1}-N_2(I-\ri\e N_2)^{-1})g(N_2)\\[.2cm]
&=(I-\ri\e N_1)^{-1}(f(N_1)-f(N_2))-(I-\ri\e N_1)^{-1}(N_1-N_2)g(N_2)\\[.2cm]
&+(I-\ri\e N_1)^{-1}(N_1-N_2)(I-\ri\e N_2)^{-1}g(N_2).
\end{align*}
Now it is clear that $f_\e\in\OLA(\C_+)$
and we have
$$
f_\e(L)-f_\e(M)=\iint_{\R\times\R}(\dg f_\e)(x,y)\,d\mE_L(x)(L-M)\,d\mE_M(y)
$$
for every $\e>0$.

We can still assume that $f(0)=0$. Then $|f(z)|\le k|z|$ everywhere in $\C_+$ for some 
$k>0$. Hence, $|f_\e(z)|\le k\e^{-1}$ for all $\e>0$ and $z\in\C_+$.

Let us observe that  $\lim\limits_{\e\to0^+}f_\e(L)x=f(L)x$ for every $x$ in $\cd_L$.
Indeed, keeping in mind that  $|(1-\ri\e z)^{-1}|\le1$ and $\lim\limits_{\e\to0^+}(1-\ri\e z)^{-1}=1$ whenever $\im z\ge0$, we see that
$$
\lim_{\e\to0^+}(I-\ri\e L)^{-1}=\lim_{\e\to0^+}\int_\R(1-\ri\e x)^{-1}\,d\mE_L(x)=I
$$ 
in the strong
operator topology. Hence, 
$$
\lim_{\e\to0^+}(I-\ri\e L)^{-1}f(L)(L+\ri I)^{-1}=f(L)(L+\ri I)^{-1}
$$
in the strong operator topology. It follows that $\lim\limits_{\e\to0^+}f_\e(L)u=f(L)u$
for all $u\in\cd_L$. Similarly, $\lim\limits_{\e\to0^+}f_\e(M)u=f(M)u$
for all $u\in\cd_M=\cd_L$.

To prove that $\lim\limits_{\e\to0^+}(f_\e(L)-f_\e(M))=f(L)-f(M)$ in the strong operator topology
it suffices to verify that $\|f_\e(L)-f_\e(M)\|\le c\|L-M\|$, where $c$ can depend only on $f$. 

We have
\begin{align}
\label{Dfe}
(\dg f_\e)(x,y)&=\frac{(1-\ri\e x)^{-1}f(x)-(1-\ri\e y)^{-1}f(y)}{x-y}\nonumber\\[.2cm]
&=\frac{f(x)-\ri\e yf(x)-f(y)+\ri\e xf(y)}{(1-\ri\e x)(x-y)(1-\ri\e y)}\nonumber\\[.2cm]
&=(\dg f)(x,y)(1-\ri\e y)^{-1}+\ri\e f_\e(x)(1-\ri\e y)^{-1}.
\end{align}
Hence, 
\begin{align*}
\|(\dg f_\e)(x,y)\|_{C_{A,\be}(\R)\otimes_{\rm h}C_{A,\be}(\R)}
&\le\|(\dg f)(x,y)(1-\ri\e y)^{-1}\|_{C_{A,\be}(\R)\otimes_{\rm h}C_{A,\be}(\R)}\\[.2cm]
&+\|\ri\e f_\e(x)(1-\ri\e y)^{-1}\|_{C_{A,\be}(\R)\otimes_{\rm h}C_{A,\be}(\R)}\\[.2cm]
&\le\|(\dg f)(x,y)\|_{C_{A,\be}(\R)\otimes_{\rm h}C_{A,\be}(\R)}+\e\|f_\e\|_{C_{A,\be}(\R)}\\[.2cm]
&\le\|(\dg f)(x,y)\|_{C_{A,\be}(\R)\otimes_{\rm h}C_{A,\be}(\R)}+k,
\end{align*}
whence
$\|f_\e(L)-f_\e(M)\|\le\left(\|(\dg f)(x,y)\|_{C_{A,\be}(\R)\otimes_{\rm h}C_{A,\be}(\R)}+k\right)\|L-M\|$.

It remains to prove that
$$
\lim_{\e\to0^+}\iint\limits_{\R\times\R}(\dg f_\e)(x,y)\,d\mE_L(x)(L-M)\,d\mE_M(y)=
\iint\limits_{\R\times\R}(\dg f)(x,y)\,d\mE_L(x)(L-M)\,d\mE_M(y)
$$
in the strong operator topology.
Applying \rf{Dfe}, we obtain
\begin{align*}
f_\e(L)&-f_\e(M)=
\iint\limits_{\R\times\R}(\dg f_\e)(x,y)\,d\mE_L(x)(L-M)\,d\mE_M(y)\\[.2cm]
&=\left(\,\,\iint\limits_{\R\times\R}(\dg f)(x,y)\,d\mE_L(x)(L-M)\,d\mE_M(y)+\ri\e f_\e(L)(L-M)\right)(1-\ri\e M)^{-1}.
\end{align*}
Now to pass to the limit as $\e\to0^+$ it suffices to observe that
$\lim\limits_{\e\to0^+}(1-\ri\e M)^{-1}=I$ and $\lim\limits_{\e\to0^+}\e f_\e(L)=\0$ in the strong operator topology. 
The first equality has been proved above while the second equality can be proved in the same way. $\bl$

\ 

\section{\bf Divided differences and Haagerup tensor products}
\setcounter{equation}{0}
\label{DdiHtp}

\

The main result of the paper will be obtained in \S\:\ref{dvepary}. To obtain such an estimate, we are going to use an integral formula for the operator difference
$f(L_1,M_1)-f(L_2,M_2)$, where $(L_1,M_1)$ and $(L_2,M_2)$ are pairs of commuting maximal dissipative operators and $f\in\big(B_{\be,1}^1\big)_+(\R^2)$. We are going to establish the following formula:
\begin{align}
\label{glafor}
f(L_1,M_1)&-f(L_2,M_2)=\iint\limits_{\R^2\times\R^2}
\frac{f(s_1,t_1)-f(s_1,t_2)}{t_1-t_2}\,d\mE_1(s_1,t_1)(M_1-M_2)\,d\mE_2(s_2,t_2)
\nonumber\\[.2cm]
&+\iint\limits_{\R^2\times\R^2}
\frac{f(s_1,t_2)-f(s_2,t_2)}{s_1-s_2}\,d\mE_1(s_1,t_1)(L_1-L_2)\,d\mE_2(s_2,t_2),
\end{align}
where $\mE_1$ is a semi-spectral measure of $(L_1,M_1)$ and $\mE_2$ is a semi-spectral measure of $(L_2,M_2)$.

It is easy to see that it suffices to prove \rf{glafor} for functions $f$ in the class
$(\E^\be_\s)_+(\R^2)$ defined below.

To establish this formula, we first have to show that the integrands in both double operator integrals on the right belong to the Haagerup integral product 
$L^\be_{\E_1}\otimes_{\rm h}L^\be_{\E_1}$, and so the right-hand side of \rf{glafor} makes sense. Then we have to prove that the right-hand side of \rf{glafor} coincides with its left-hand side.

To this end we show that
\begin{align*}
\iint\limits_{\R^2\times\R^2}
\frac{f(s_1,t_1)-f(s_1,t_2)}{t_1-t_2}\,d\mE_1(s_1,t_1)(M_1-M_2)\,d\mE_2(s_2,t_2)
\\[.2cm]
=f(L_1,M_1)-\iint\limits_{\R\times\R}f(s,t)\,d\mE_{L_1}(s)\,d\mE_{M_2}(t)
\end{align*}
and
\begin{align*}
\iint\limits_{\R^2\times\R^2}
\frac{f(s_1,t_2)-f(s_2,t_2)}{s_1-s_2}\,d\mE_1(s_1,t_1)(L_1-L_2)\,d\mE_2(s_2,t_2)
\\[.2cm]
=\iint\limits_{\R\times\R}f(s,t)\,d\mE_{L_1}(s)\,d\mE_{M_2}(t)-f(L_2,M_2).
\end{align*}
Note that both identities involve the double operator integral
$$
\iint\limits_{\R\times\R}f(s,t)\,d\mE_{L_1}(s)\,d\mE_{M_2}(t).
$$
It can be interpreted as a function $f(L_1,M_2)$ of our ({\it not necessarily commuting}) maximal dissipative operators $L_1$ and $M_2$. We define in \S\:\ref{nekommu} such functions of noncommutative maximal dissipative operators and we will prove the above identities in \S\:\ref{dvepary}.

Let $\s\in(0,\be)$ and $p\in[1,\be]$. Denote by $\E^p_\s(\R^d)$ the space of the functions $f\in L^p(\R^d)$
such that $\supp\F f\subset\s B^d$, where $B^d$ denotes the closed unit ball of $\R^d$.  It is well known that each function $f\in\E^p_\s(\R^d)$
is the restriction to $\R^d$ of an entire function on $\C^d$ of exponential type at most $\s$. Put
$$
(\E^p_\s)_+(\R^d)=\{f\in L^p(\R^d):~\supp\F f\subset[0,\be)^d\}.
$$
We are going to consider the cases $d=1$ and $d=2$. In the case $d=1$ we write $\E^p_\s$ and $(\E^p_\s)_+$
instead $\E^p_\s(\R)$ and $(\E^p_\s)_+(\R)$.

%
%
%

It follows from the results of \cite{Pe3} that
\bay
\label{m1}
f\in\E_\s^\be\quad\Longrightarrow\quad
\left\|\frac{f(x)-f(y)}{x-y}\right\|_{\frak M(E_1,E_2)}\le\const\s\|f\|_{L^\be(\R)}
\ey
for all Borel spectral measures $E_1$ and $E_2$ on $\R$.
It was shown in \cite{AP4} that inequality \rf{m1} holds with constant equal to 1.

In \cite{APPS} for $f\in\E^\be_\s$, an explicit representation of the divided difference $\frac{f(x)-f(y)}{x-y}$  was obtained
as an element of ${\rm C}_{\rm b}(\R)\otimes_{\rm h}{\rm C}_{\rm b}(\R)$,
where ${\rm C}_{\rm b}(\R)\df L^\infty(\R)\cap {\rm C}(\R)$.
Clearly such a representation holds for all functions $f$ in
 $(\E_\s^\be)_+$.

However, we need an expansion $\frac{f(x)-f(y)}{x-y}=\sum_n\f_n(x)\psi_n(y)$ with the 
additional  requirement that $\f_n$ and $\psi_n\in H^\be(\C_+)$.

\begin{thm}
\label{to}
Let $f\in(\E_\s^\be)_+$. Then
\begin{align}
\label{haa1}
\frac{f(x)-f(y)}{x-y}&=\sum_{n\in\Z}\s\cdot\frac{f(x)-f(2\s^{-1}\pi n)}{\s x-2\pi n}\cdot\frac{e^{{\rm i}\s y}-1}{{\rm i}(\s y-2\pi n)}\\[.2cm]
\label{haa2}
&=\frac1{2\pi\rm i}\int_\R\frac{f(x)-f(t)}{x-t}\cdot\frac{e^{{\rm i}\s(y-t)}-1}{y-t}\,dt,
\quad x,~y\in\R.
\end{align}
Moreover,
\bay
\label{2vy}
\sum_{n\in\Z}\frac{\big|f(x)-f\big(2\pi\s^{-1}n\big)\big|^2}{(\s x-2\pi n)^2}
=\frac1{2\pi\s}\int_\R\frac{|f(x)-f(t)|^2}{(x-t)^2}dt\le\frac4\pi\|f\|_{L^\be(\R)}^2,
\quad x\in\R,
\ey
and
\bay
\label{hiz}
\sum_{n\in\Z}\frac{|e^{{\rm i}\s y}-1|^2}{(\s y-2\pi n)^2}=1=\frac1{2\pi\s}\int_\R\frac{|e^{{\rm i}\s(y-t)}-1|^2}{(y-t)^2}\,dt,\quad y\in\R.
\ey
\end{thm}

Note that the identities in \rf{hiz} are elementary and well known.

\medskip

\Pf Clearly, it suffices to consider the case $\s=1$. 

We are going to use the well-known facts that the family
$\left\{\frac{e^{{\rm i}z}-1}{{\rm i}(z-2\pi n)}\right\}_{n\in\Z}$
forms an orthogonal basis in the space $(\E_1^2)_+$,
\bay
\label{Fz}
g(z)=\sum_{n\in\Z}g(2\pi n)\dfrac{e^{{\rm i}z}-1}{{\rm i}(z-2\pi n)},
\ey
and
\bay
\label{Fz1}
\sum_{n\in\Z}|g(2\pi n)|^2=\frac1{2\pi}\int_\R|g(t)|^2\,dt.
\ey
for every $g\in(\E_1^2)_+$  (this follows from Th. 1 in \cite{L}, Lect. 20.2). It follows immediately from \rf{Fz1} that
\bay
\label{Fz2}
\sum_{n\in\Z}g(2\pi n)\ov{h(2\pi n)}=\frac1{2\pi}\int_\R g(t)\ov{h(t)}\,dt
\quad\mbox{for every}\quad g,~h\in(\E_1^2)_+.
\ey

Given $x\in\R$, we consider the function $g$ defined by $g(\l)=\frac{f(x)-f(\l)}{x-\l}$, $\l\in\C$. Clearly, $g\in(\E_1^2)_+$.

It is easy to see that \rf{haa1} is a consequence of
\rf{Fz} and the equality in \rf{2vy} is a consequence of \rf{Fz1}.
It is also easy to see that \rf{haa2} follows from \rf{Fz2}.

It remains to prove that
$$
\frac1{\pi}\int_\R\frac{|f(x)-f(t)|^2}{(x-t)^2}dt\le\frac4\pi\|f\|_{L^\be(\R)}^2
$$
for arbitrary $f\in(\E_1^\be)_+$ and $x\in\R$. Without loss of generality we may assume that
$\|f\|_{L^\be(\R)}=1$. Then $\|f^\prime\|_{L^\be(\R)}\le1$ by the Bernstein
inequality. Hence, $|f(x)-f(t)|\le\min\{2,|x-t|\}$, and we have
\bey
\frac1{2\pi}\int_\R\frac{|f(x)-f(t)|^2}{(x-t)^2}dt\le\frac1{2\pi}\int_\R\frac{\min\{4,(x-t)^2\}}{(x-t)^2}dt
=\frac1\pi\int_0^2dt+\frac4\pi\int_2^\be\frac{dt}{t^2}=\frac4\pi.\quad\bl
\eey

\medskip

{\bf Remark.} Note that the equality
$$
\frac{f(x)-f(y)}{x-y}=\frac1{2\pi\rm i}\int_\R\frac{f(x)-f(t)}{x-t}\cdot\frac{e^{{\rm i}\s(y-t)}-1}{y-t}\,dt,
$$
is an immediate consequence of the well-known fact that
$\dfrac{e^{{\rm i}\s(x-y)}-1}{2\pi{\rm i}(x-y)}$ is the reproducing kernel for the
functional Hilbert space $(\E_\s^2)_+$.

\medskip


Let $f\in(\E_\s^\be)_+(\R^2)$.
Put 
$$
\big(\dg_xf\big)(x_1,y_1;x_2,y_2)\df\big(\dg_xf\big)(x_1;x_2,y_2)\df\frac{f(x_1,y_2)-f(x_2,y_2)}{x_1-x_2}
$$
and
$$
\big(\dg_yf\big)(x_1,y_1;x_2,y_2)\df\big(\dg_yf\big)(x_1,y_1;y_2)\df\frac{f(x_1,y_1)-f(x_1,y_2)}{y_1-y_2}.
$$

Let $\CAbe(\R^d)\df\{f\in{\rm C}(\R^d)\cap L^\be(\R^d): \supp\F f\subset[0,\be)^d\}$, $\CAbe\df \CAbe(\R)=\CAbe(\C_+)$.

\begin{thm}
\label{fs}
Let $\s>0$ and let $f\in(\E_\s^\be)_+(\R^2)$. 
Then
$$
\dg_xf,\,\dg_yf\in \CAbe(\R^2)\otimes_{\rm h}\CAbe(\R^2),
$$
$$
\|\dg_xf\|_{\CAbe(\R^2)\otimes_{\rm h}\CAbe(\R^2)}\le\frac2{\sqrt\pi}\,\s\|f\|_{L^\be(\R^2)}
$$
and
$$
\|\dg_yf\|_{\CAbe(\R^2)\otimes_{\rm h}\CAbe(\R^2)}\le\frac2{\sqrt\pi}\,\s\|f\|_{L^\be(\R^2)}.
$$
\end{thm}

\Pf Clearly, $f$ is the restriction to $\R^2$ of an entire function of two complex variables.
Moreover, $f(\cdot,a),\,f(a,\cdot)\in \E_\s^\be$ for every $a\in\R$.
It suffices to consider the case $\s=1$. By Theorem \ref{to}, we have
\bey
\big(\dg_xf\big)(x_1,y_1;x_2,y_2)=\frac{f(x_1,y_2)-f(x_2,y_2)}{x_1-x_2}=
\sum_{n\in\Z}\frac{e^{{\rm i}x_1}-1}{{\rm i}(x_1-2\pi n)}\cdot\frac{f(x_2,y_2)-f(2\pi n,y_2)}{x_2-2\pi n}
\eey
and

\bey
\big(\dg_yf\big)(x_1,y_1;x_2,y_2)=\frac{f(x_1,y_1)-f(x_1,y_2)}{y_1-y_2}=\sum_{n\in\Z}\frac{f(x_1,y_1)-f(x_1,2\pi n)}{y_1-2\pi n}
\cdot\frac{e^{{\rm i}y_2}-1}{{\rm i}(y_2-2\pi n)}.
\eey
Moreover, by Theorem \ref{to},
we have
\bey
\sum_{n\in\Z}\frac{|f(x_1,y_1)-f(x_1,2\pi n)|^2}{(y_1-2\pi n)^2}\le\frac4\pi\|f(x_1,\cdot)\|_{L^\be(\R)}^2
\le3\|f\|_{L^\be(\R^2)}^2,\\
\sum_{n\in\Z}\frac{|f(x_2,y_2)-f(2\pi n,y_2)|^2}{(x_2-2\pi n)^2}\le\frac4\pi\|f(\cdot,y_2)\|_{L^\be(\R)}^2
\le3\|f\|_{L^\be(\R^2)}^2,\\
\eey
and
$$
\sum_{n\in\Z}\frac{|e^{{\rm i}x_1}-1|^2}{(x_1-2\pi n)^2}=\sum_{n\in\Z}\frac{|e^{{\rm i}y_2}-1|^2}{(y_2-2\pi n)^2}=1.
$$
This implies the result. $\bl$

%
%
%
%

\medskip

Put $\vk(x)\df(1-\ri x)^{-1}$.  

\begin{cor} 
\label{cor23}
Let $f\in(\E_\s^\be)_+(\R^2)$, where $\s>0$. Then 
$$
\vk(x)f(x,y)\in  \CAbe\otimes_{\rm h}\CAbe
\quad\text{and}\quad \vk(y)f(x,y)\in \CAbe\otimes_{\rm h}\CAbe.
$$
Moreover,
$$
\|\vk(x)f(x,y)\|_{\CAbe\otimes_{\rm h}\CAbe}\le\const(1+\s)\|f\|_{L^\be(\R^2)}
$$
and
$$
\|\vk(y)f(x,y)\|_{\CAbe\otimes_{\rm h}\CAbe}\le\const(1+\s)\|f\|_{L^\be(\R^2)}.
$$
\end{cor}

\Pf It suffices to prove the result for the function $\vk(x)f(x,y)$. 
We have $f(x,y)=x\big(\dg_xf\big)(x;0,y)+f(0,y)$. Hence, 
\bey
\|\vk(x)f(x,y)\|_{\CAbe\otimes_{\rm h}\CAbe}\le
\|(1-\ri x)^{-1}x\big(\dg_xf\big)(x;0,y)\|_{\CAbe\otimes_{\rm h}\CAbe}\\
+\|(1-\ri x)^{-1}f(0,y)\|_{\CAbe\otimes_{\rm h}\CAbe}\le\frac2{\sqrt\pi}\s\|f\|_{L^\be(\R^2)}+
\|f\|_{L^\be(\R^2)}.\quad\bl
\eey

%

\

\section{\bf Functions of pairs of noncommuting maximal dissipative operators}
\setcounter{equation}{0}
\label{nekommu}

\

In this section we consider functions of noncommuting maximal dissipative operators.
Let $(L,M)$ be a pair of not necessarily commuting 
maximal dissipative operators in a Hilbert space $\h$.

For $f\in \CAbe\otimes_{\rm h}\CAbe$,
we put
$$
f(L,M)\df\iint\limits_{\R\times\R}f(s,t)\,d\mE_L(s)\,d\mE_M(t).
$$

Clearly, 
$\|f(L,M)\|\le\|f\|_{\CAbe\otimes_{\rm h}\CAbe}$.
It is easy to see that if $f\in \CAbe\otimes_{\rm h}\CAbe$, $g\in \CAbe$ and $f_\sharp$ is the function defined by
$f_\sharp(s,t)\df g(t)f(s,t)$, then $f_\sharp\in\CAbe\otimes_{\rm h}\CAbe$ and $f_\sharp(L,M)=f(L,M)g(M)$.
This gives us an idea how to define $f(L,M)$ in a more general case.

Suppose that $f_\sharp(s,t)\df(1-\ri t)^{-1}f(s,t)\in \CAbe\otimes_{\rm h}\CAbe$. Put
$$
f(L,M)\df f_\sharp(L,M)(I-\ri M)=\left(\,\,\,\iint\limits_{\R\times\R}f_\sharp(s,t)\,d\mE_L(s)\,d\mE_M(t)\right)(I-\ri M).
$$
Note that $f(L,M)$  is not necessarily a bounded operator but the operator $f(L,M)\lb(I-\ri M)^{-1}$ is bounded.

\medskip

{\bf Remark.} Corollary \ref{cor23} implies now that $f(L,M)$ is defined for every $f\in(\E_\s^\be)_+(\R^2)$, where $\s>0$.
Moreover, $f(L,M)(I-\ri M)^{-1}$ is bounded for such functions.

\

\section{\bf The principal inequality}
\setcounter{equation}{0}
\label{dvepary}

\

In this section we obtain the main results of the paper. 

\begin{lem}
\label{l1m2}
Let $f\in(\E_\s^\be)_+(\R^2)$, where $\s>0$. Suppose that $(L_1,M_1)$ and $(L_2,M_2)$ are pairs of commuting 
maximal dissipative operators such that the operators $L_1-L_2$ and $M_1-M_2$ are bounded. Then the following formulas hold
\begin{align}
\label{31}
f(L_1,M_1)&-f(L_1,M_2)\nonumber\\[.2cm]
&=\iint\limits_{\R^2\times\R^2}
\frac{f(s_1,t_1)-f(s_1,t_2)}{t_1-t_2}\,d\mE_1(s_1,t_1)(M_1-M_2)\,d\mE_2(s_2,t_2)
\end{align}
and
\begin{align}
\label{32}
f(L_1,M_2)&-f(L_2,M_2)\nonumber\\[.2cm]
&=\iint\limits_{\R^2\times\R^2}
\frac{f(s_1,t_2)-f(s_2,t_2)}{s_1-s_2}\,d\mE_1(s_1,t_1)(L_1-L_2)\,d\mE_2(s_2,t_2),
\end{align}
where $\mE_1$ and $\mE_2$ are semi-spectral measures of the pairs
$(L_1,M_1)$ and $(L_2,M_2)$.

Moreover, the operator $f(L_1,M_2)$ is bounded.
\end{lem}

\Pf First, we prove equality \rf{31}. Let $X$ and $Y$ be the left-hand side and the right-hand side of  \rf{31}.
Clearly, $Y$ is a bounded operator. Note that $f(L_1,M_1)$ is bounded (see \S\ref{Disopfuis}).
To prove that $f(L_1,M_2)$ is bounded, we first recall that $f(L_1,M_2)(I-\ri M_2)^{-1}$ is bounded, see \S\:\ref{nekommu}.
Clearly, $Y(I-\ri M_2)^{-1}$ is a bounded operator. To prove that $X=Y$
it suffices to verify that $(I-\ri M_1)^{-1}(X(I-\ri M_2)^{-1})=(I-\ri M_1)^{-1}(Y(I-\ri M_2)^{-1})$.
We have 
\begin{multline*}
(I-\ri M_1)^{-1}f(L_1,M_1)(I-\ri M_2)^{-1}-(I-\ri M_1)^{-1}f(L_1,M_2)(I-\ri M_2)^{-1}\\
=\iint\limits_{\R^2\times\R^2}(f(s_1,t_1)-f(s_1,t_2))(1-\ri t_1)^{-1}(1-\ri t_2)^{-1}\,d\mE_1(s_1,t_1)\,d\mE_2(s_2,t_2)\\
=\iint\limits_{\R^2\times\R^2}\frac{f(s_1,t_1)-f(s_1,t_2)}{t_1-t_2}
\left(\frac{t_1}{1-\ri t_1}-\frac{t_2}{1-\ri t_2}\right)\,d\mE_1(s_1,t_1)\,d\mE_2(s_2,t_2)\\
=\iint\limits_{\R^2\times\R^2}\frac{t_1}{1-\ri t_1}\cdot\frac{f(s_1,t_1)-f(s_1,t_2)}{t_1-t_2}\,d\mE_1(s_1,t_1)\,d\mE_2(s_2,t_2)\\
-\iint\limits_{\R^2\times\R^2}\frac{f(s_1,t_1)-f(s_1,t_2)}{t_1-t_2}\cdot\frac{t_2}{1-\ri t_2}\,d\mE_1(s_1,t_1)\,d\mE_2(s_2,t_2)\\
=\iint\limits_{\R^2\times\R^2}\frac{f(s_1,t_1)-f(s_1,t_2)}{t_1-t_2}\,d\mE_1(s_1,t_1)M_1(I-\ri M_1)^{-1}\,d\mE_2(s_2,t_2)\\
-\iint\limits_{\R^2\times\R^2}\frac{f(s_1,t_1)-f(s_1,t_2)}{t_1-t_2}\,d\mE_1(s_1,t_1)M_2(I-\ri M_2)^{-1}\,d\mE_2(s_2,t_2)\\
=(I-\ri M_1)^{-1}\left(\,\,\,\iint\limits_{\R^2\times\R^2}\frac{f(s_1,t_1)-f(s_1,t_2)}{t_1-t_2}\,d\mE_1(s_1,t_1)(M_1-M_2)\,d\mE_2(s_2,t_2)\right)
(I-\ri M_2)^{-1}
\end{multline*}
which implies \rf{31}.
In the same way we can verify that
\bey
(I-\ri L_1)^{-1}f(L_1,M_1)(I-\ri L_2)^{-1}-(I-\ri L_1)^{-1}f(L_1,M_2)(I-\ri L_2)^{-1}\\
=(I-\ri L_1)^{-1}\left(\,\,\,\iint\limits_{\R^2\times\R^2}\frac{f(s_1,t_2)-f(s_2,t_2)}{s_1-s_2}\,d\mE_1(s_1,t_1)(L_1-L_2)\,d\mE_2(s_2,t_2)\right)
(I-\ri L_2)^{-1},
\eey
which implies \rf{32}.   

To see that the operator $f(L_1,M_2)$ is bounded, it suffices to observe that 
 $f(L_1,M_1)$ is bounded and the double operator integral on the right-hand side of \rf{31} is also bounded. $\bl$

\begin{cor} 
Under the hypotheses
of Lemma {\rm\ref{l1m2}}, the following formula holds:
\begin{align*}
f(L_1,M_1)-f(L_2,M_2)&=\iint\limits_{\R^2\times\R^2}
\frac{f(s_1,t_1)-f(s_1,t_2)}{t_1-t_2}\,d\mE_1(s_1,t_1)(M_1-M_2)\,d\mE_2(s_2,t_2)\\[.2cm]
&+\iint\limits_{\R^2\times\R^2}
\frac{f(s_1,t_2)-f(s_2,t_2)}{s_1-s_2}\,d\mE_1(s_1,t_1)(L_1-L_2)\,d\mE_2(s_2,t_2).
\end{align*}
\end{cor}

\begin{cor} 
\label{43}
Under the hypotheses of Lemma {\rm\ref{l1m2}},
$$
\|f(L_1,M_1)-f(L_2,M_2)\|\le\const\s\|f\|_{L^\be(\R^2)}\max(\|L_1-L_2\|,\|M_1-M_2\|).
$$
\end{cor}

Let $f\in\big(B_{\be,1}^1\big)_+(\R^2)$. Then 
$f(\z_1,\z_2)(\z_1+\ri)^{-1}(\z_2+\ri)^{-1}\in\CAb$, and $f(L,M)$ is defined for
any pair $(L,M)$ of commuting maximal dissipative operators, see \S\ref{In}.

The following theorem is the main result of the paper.

\begin{thm}
\label{osnrez}
Let $f\in\big(B_{\be,1}^1\big)_+(\R^2)$. Suppose that $(L_1,M_1)$ and $(L_2,M_2)$ are pairs of commuting 
maximal dissipative operators such that the operators $L_1-L_2$ and $M_1-M_2$ are bounded. Then 
$$
\|f(L_1,M_1)-f(L_2,M_2)\|\le\const\|f\|_{(B_{\be,1}^1)(\R^2)}\max\{\|L_1-L_2\|,\|M_1-M_2\|\}.
$$
\end{thm}

\Pf Consider the functions $f_n$, $n\ge0$, defined by \rf{fn}. Clearly, by Corollary \ref{43},
\begin{align*}
\|f(L_1,M_1)-f(L_2,M_2)\|&\le\sum_{n\ge0}\|f_n(L_1,M_1)-f_n(L_2,M_2)\|\\
&\le\const\left(\sum_{n\ge0}2^n\|f_n\|_{L^\be(\R^2)}\right)\max\{\|L_1-L_2\|,\|M_1-M_2\|\}
\\[.2cm]
&\le\const\|f\|_{(B_{\be,1}^1)(\R^2)}\max\{\|L_1-L_2\|,\|M_1-M_2\|\}.\quad\bl
\end{align*}

\

\section{\bf H\"older type estimates and estimates in Schatten--von Neumann norms}
\setcounter{equation}{0}
\label{HolSp}

\

In this section we use the main result of the previous section to obtain  H\"older type estimates and estimates in Schatten--von Neumann norms for functions of pairs of commuting dissipative operators under perturbation.

Let $\a\in(0,1)$. Recall that the {\it analytic H\"older class} $(\L_\a(\C_+))_+$ 
is defined as the space of functions $f$ analytic in $\C_+$, continuous on $\clos\C_+$ and satisfying the inequality
$$
|f(\z)-f(\t)|\le\const|\z-\t|^\a,\quad\z,~\t\in\clos\C_+.
$$

\begin{thm}
\label{Gelder}
Let $\a\in(0,1)$ and let $f\in(\L_\a(\C_+))_+$. Suppose that $(L_1,M_1)$ and $(L_2,M_2)$ are pairs of commuting 
maximal dissipative operators such that the operators $L_1-L_2$ and $M_1-M_2$ are bounded.  Then 
$$
\|f(L_1,M_1)-f(L_2,M_2)\|\le\const\|f\|_{(\L_\a(\C_+))_+}
\max\{\|L_1-L_2\|^\a,\|M_1-M_2\|^\a\}.
$$
\end{thm}

Let us consider now a more general case. Let $\o$ be a {\it modulus of continuity}, i.e., $\o$ is a real continuous nondecreasing function on $[0,\be)$ such that
$\o(t+s)\le\o(t)+\o(s)$, $s,\,t\ge0$. The space $(\L_\o(\C_+))_+$ is, by definition, the class of functions $f$ analytic in $\C_+$, continuous on $\clos\C_+$ and satisfying the inequality
$$
|f(\z)-f(\t)|\le\const\o(|\z-\t|),\quad\z,~\t\in\clos\C_+.
$$
Given a modulus of continuity $\o$, we put
$$
\o_*(s)\df s\int_s^\be\frac{\o(t)}{t^2}\,dt=\int_1^\be\frac{\o(st)}{t^2}\,dt,\quad s>0.
$$
It is easy to see that if $\o_*(s)<\be$ for some $s>0$, then $\o_*(s)<\be$ for all $s>0$ and
$\o_*$ is a modulus of continuity.

\begin{thm}
\label{Lomega}
Let $\o$ be a modulus of continuity and let $f\in(\L_\o(\C_+))_+$. Suppose that $(L_1,M_1)$ and $(L_2,M_2)$ are pairs of commuting 
maximal dissipative operators such that the operators $L_1-L_2$ and $M_1-M_2$ are bounded.  Then 
$$
\|f(L_1,M_1)-f(L_2,M_2)\|\le\const\|f\|_{(\L_\o(\C_+))_+}
\max\o_*\big(\{\|L_1-L_2\|,\|M_1-M_2\|\}\big).
$$
\end{thm}

The next theorem allows us to obtain estimates in Schatten--von Neumann norms $\bS_p$.

\begin{thm}
\label{SpGeld}
Let $0<\a<1$ and  $p>1$. Suppose that $(L_1,M_1)$ and $(L_2,M_2)$ are pairs of commuting maximal dissipative operators such that  $L_1-L_2\in\bS_p$ and $M_1-M_2\in\bS_p$. Then $f(L_1,M_1)-f(L_2,M_2)\in\bS_{p/\a}$ and
$$
\|f(L_1,M_1)-f(L_2,M_2)\|_{\bS_{p/\a}}\le\const\|f\|_{(\L_\a(\C_+))_+}
\max\{\|L_1-L_2\|_{\bS_p}^\a,\|M_1-M_2\|_{\bS_p}^\a\}.
$$
\end{thm}

Theorems \ref{Gelder}, \ref{Lomega} and \ref{SpGeld} can be deduced from Theorem \ref{osnrez} in the same way as it was done for functions of self-adjoint operators, see
\cite{AP1} and \cite{AP2}.

\

\
 
 \begin{footnotesize}
 
\noindent
\begin{tabular}{p{7cm}p{15cm}}
A.B. Aleksandrov & V.V. Peller \\
St.Petersburg Branch & Department of Mathematics\\
Steklov Institute of Mathematics  & and Computer Sciences\\
Fontanka 27, 191023 St.Petersburg & St.Petersburg State University\\
Russia&Universitetskaya nab., 7/9,\\
email: alex@pdmi.ras.ru&199034 St.Petersburg, Russia\\
&\\
&Department of Mathematics\\
&Michigan State University\\
&East Lansing, Michigan 48824\\
&USA\\
&and\\
&Peoples' Friendship University\\
& of Russia (RUDN University)\\
&6 Miklukho-Maklaya St., Moscow,\\
& 117198, Russian Federation\\
& email: peller@math.msu.edu
\end{tabular}

\end{footnotesize}

\end{document}

%% file: classstartscr.tex
\newcommand{\lb}{\linebreak}

\renewcommand{\a}{\alpha}

\newcommand{\e}{\varepsilon}
\newcommand{\vk}{\varkappa}
\newcommand{\z}{\zeta}

\renewcommand{\l}{\lambda}

\newcommand{\s}{\sigma}
\renewcommand{\t}{\tau}

\newcommand{\f}{\varphi}
\renewcommand{\o}{\omega}

\newcommand{\D}{\Delta}
\renewcommand{\L}{\Lambda}

\newcommand{\E}{{\mathscr E}}
\newcommand{\cd}{{\mathscr D}}
\newcommand{\F}{{\mathscr F}}

\newcommand{\h}{{\mathscr H}}

\newcommand{\K}{{\mathscr K}}
\newcommand{\cL}{{\mathscr L}}
\newcommand{\M}{{\mathscr M}}
\newcommand{\N}{{\mathscr N}}

\newcommand{\X}{{\mathscr X}}

\newcommand{\C}{{\Bbb C}}
\newcommand{\T}{{\Bbb T}}

\newcommand{\dd}{{\Bbb D}}
\newcommand{\R}{{\Bbb R}}
\newcommand{\Z}{{\Bbb Z}}

\newcommand{\0}{{\boldsymbol{0}}}

\newcommand{\bs}{\boldsymbol}

\newcommand{\bS}{{\boldsymbol S}}

\newcommand{\rf}[1]{(\ref{#1})}

\newcommand{\df}{\stackrel{\mathrm{def}}{=}}

\newcommand{\Ker}{\operatorname{Ker}}

\newcommand{\spn}{\operatorname{span}}
\newcommand{\supp}{\operatorname{supp}}
\newcommand{\clos}{\operatorname{clos}}

\newcommand{\const}{\operatorname{const}}

\newcommand{\eeq}{\end{equation}}
\newcommand{\beq}{\begin{equation}}
\newcommand{\bay}{\begin{eqnarray}}
\newcommand{\ba}{\begin{align*}}
\newcommand{\ea}{\end{align*}}
\newcommand{\ey}{\end{eqnarray}}
\newcommand{\bey}{\begin{eqnarray*}}
\newcommand{\eey}{\end{eqnarray*}}

\newcommand{\be}{\infty}

\newcommand{\bl}{\blacksquare}

\newcommand{\Range}{\operatorname{Range}}

\newcommand{\Pf}{{\bf Proof. }}
\newcommand{\im}{\operatorname{Im}}

\newcommand{\ov}{\overline}

\newtheorem{thm}{\hspace{\parindent}Theorem}[section]

\newtheorem{cor}[thm]{\hspace{\parindent}Corollary}
\newtheorem{lem}[thm]{\hspace{\parindent}Lemma}

%% file: kommutir_dissipativ03.bbl
\begin{thebibliography}{99}

\bibitem[A]{A}
{\sc A. B. Aleksandrov}, {\em Operator Lipschitz functions and linear fractional transformations}, Zap. Nauchn. Sem. 
POMI {\bf401} (2012), Issledovaniya po Lineinym Operatoram i Teorii Funktsii. {\bf40}, 5--52;

English transl.: J. Math. Sci. (N. Y.), {\bf194:6} (2013),  603--627.

\bibitem[ANP]{ANP} {\sc A.B. Aleksandrov, F.L. Nazarov} and
{\sc V.V. Peller}, {\em Functions of noncommuting self-adjoint operators under perturbation and estimates of triple operator integrals},  Adv. Math. {\bf295} (2016), 1--52.

\bibitem[AP1]{AP1}  {\sc A.B. Aleksandrov} and {\sc V.V. Peller},  {\em Operator H\"older--Zygmund functions}, Advances in Math.
{\bf224} (2010), 910--966.

\bibitem[AP2]{AP2}  {\sc A.B. Aleksandrov} and {\sc V.V. Peller},  {\em Functions of operators under perturbations of
class $\bS_p$}, J. Funct. Anal. {\bf258} (2010), 3675--3724.

\bibitem[AP3]{AP3}  {\sc A.B. Aleksandrov} and {\sc V.V. Peller}, {\it Functions of perturbed dissipatibe operators}, Algebra i Analiz {\bf21:2} (2011),  9--51.

\noindent
English transl.: St. Petersburg Math. J. {\bf23:2} (2012),  209--238.


\bibitem[AP4]{AP4} {\sc A.B. Aleksandrov} and {\sc V.V. Peller}, {\it Operator Lipschitz functions}, 
Uspekhi Matem. Nauk {\bf71:4} (2016), 3--106.

English transl.: Russian Mathematical Surveys {\bf71:4} (2016), 605--702.

\bibitem[AP5]{AP5} {\sc A.B. Aleksandrov} and {\sc V.V. Peller}, {\it Dissipative operators and operator Lipschitz functions},
Proc. Amer. Math. Soc. {\bf147:5} (2019), 2081 -- 2093.

\bibitem[AP6]{AP6} {\sc A.B. Aleksandrov} and {\sc V.V. Peller}, {\it Functions of perturbed pairs of noncommuting contractions}, Izvestia RAN, 2020.


\bibitem[APPS]{APPS} {\sc A.B. Aleksandrov, V.V. Peller, D. Potapov}, and
{\sc F. Sukochev}, {\em Functions of normal operators under perturbations},
Advances in Math. {\bf226} (2011), 5216--5251.

\bibitem[An]{An}{\sc And\^o}, {\em On a pair of commutative contractions}, 
Acta Sci. Math. (Szeged) {\bf24} (1963), 88-90.

\bibitem[BS1]{BS1} {\sc M.S. Birman} and {\sc M.Z. Solomyak},
{\em Double Stieltjes operator integrals},
Problems of Math. Phys., Leningrad. Univ. {\bf1} (1966), 33--67 (Russian).
\newline
English transl.: Topics Math. Physics {\bf1} (1967), 25--54, Consultants Bureau Plenum
Publishing Corporation, New York.

\bibitem[BS2]{BS2} {\sc M.S. Birman} and {\sc M.Z. Solomyak},
 {\em Double Stieltjes operator integrals. II},
 Problems of Math. Phys., Leningrad. Univ. {\bf2} (1967), 26--60 (Russian).
 \newline
English transl.: Topics Math. Physics {\bf2} (1968), 19--46, Consultants Bureau Plenum
Publishing Corporation, New York.

\bibitem[BS3]{BS3} {\sc M.S. Birman} and {\sc M.Z. Solomyak},
{\em Double Stieltjes operator integrals. III},
Problems of Math. Phys., Leningrad. Univ. {\bf6} (1973), 27--53 (Russian).

\bibitem[DK]{DK} {\sc Yu.L. Daletskii} and {\sc S.G. Krein}, {\em Integration and differentiation of
functions of Hermitian operators and application to the theory of perturbations} (Russian), Trudy Sem.
Functsion. Anal., Voronezh. Gos. Univ. {\bf1} (1956), 81--105.


\bibitem[GK]{GK} {\sc I. C. Gohberg} anf {\sc M.G. Krein}, {\em Introduction to the theory of linear non-selfadjoint operators in Hilbert space}, Izdat. "Nauka'', Moscow, 1965.

\bibitem[KPSS]{KPSS} {\sc E. Kissin, D. Potapov, V. S. Shulman} and {\sc F. Sukochev}, 
{\it Operator smoothness in Schatten norms for functions of several variables: 
Lipschitz conditions, differentiability and unbounded derivations},  
Proc. Lond. Math. Soc. (3) {\bf105}:4 (2012), 661--702.

\bibitem[L]{L} {\sc B.Ya. Levin}, {\em Lectures on entire functions},
Transl. Math. Monogr., vol. 150, Amer.
Math. Soc., Providence, RI 1996, xvi+248 pp.

\bibitem[MM]{MM} {\sc M.M. Malamud} and {\sc S.M. Malamud}, {\em 
Spectral theory of operator measures in a Hilbert space},  Algebra i Analiz  {\bf15:3}  (2003), 1--77.

\bibitem[MN]{MN}
{\sc M. Malamud} and {\sc H. Neidhardt},
{\it Trace formulas for additive and non-additive perturbations},
Adv. Math. {\bf274} (2015), 736--832.

\bibitem[MNP1]{MNP1}
{\sc M.M. Malamud, H. Neidhardt} and {\sc V.V. Peller},
{\em Analytic operator Lipschitz functions in the disk and a trace formula for functions of contractions} (Russian),
Funkts. Anal. i Pril., {\bf51:3} (2017), 33--55.

English Transl.: Funct. Anal. and its Appl. {\bf51:3} (2017), 185--203.

\bibitem[MNP2]{MNP2}
{\sc M.M. Malamud, H. Neidhardt} and {\sc V.V. Peller}, {\em Absolute continuity of spectral shift}, J. Funct. Anal. {\bf276} (2019), 1575--1621.

\bibitem[Na]{Na} {\sc M.A. Naimark}, {\it Spectral functions of a symmetric operator}, 
(Russian) Bull. Acad. Sci. URSS. Ser. Math. [Izvestia Akad. Nauk SSSR] 4, (1940). 277--318.

%
%

\bibitem[Pee]{Pee} {\sc J. Peetre},
{\em New thoughts on Besov spaces}, Duke Univ. Press., Durham, NC, 1976.

\bibitem[Pe1]{Pe1} {\sc V.V. Peller},
{\em Hankel operators in the theory of perturbations of unitary and self-adjoint operators},
Funktsional. Anal. i Prilozhen. {\bf19:2}  (1985),
37--51 (Russian).

English transl.: Funct. Anal. Appl. {\bf19} (1985), 111--123.

\bibitem[Pe2]{Pe2} {\sc V.V. Peller}, {\em For which $f$ does $A-B\in{\bf S}_{p}$ 
imply that $f(A)-f(B)\in{\bf S}_{p}$?}, Operator Theory, Birkh\"{a}user,
{\bf 24} (1987), 289-294.

\bibitem[Pe3]{Pe3} {\sc V.V. Peller},
{\em Hankel operators in the perturbation theory of unbounded self-adjoint operators},
Analysis and partial differential equations,  529--544,
Lecture Notes in Pure and Appl. Math., {\bf122}, Marcel Dekker, New York, 1990.

\bibitem[Pe4]{Pe4} {\sc V.V. Peller}, {\em Differentiability of functions of contractions}, 
In: Linear and complex analysis, AMS Translations, Ser. 2 {\bf226}
(2009), 109--131, AMS, Providence.

\bibitem[Pe5]{Pe5} {\sc V.V. Peller}, {\it Multiple operator integrals in perturbation theory},  Bull. Math. Sci. {\bf6} (2016), 15--88.

\bibitem[Pe6]{Pe6} {\sc V.V. Peller}, {Functions of commuting contractions under perturbation},
Math. Nachr. {\bf291} (2019), 1151--1160.

\bibitem[PS]{PS} {\sc D. Potapov} and {\sc F. Sukochev}, {\em Operator-Lipschitz functions in Schatten--von Neumann classes}, Acta Math. {\bf207} (2011), 375--389.

\bibitem[So] {So}  {\sc B.M. Solomyak}, {\em A functional model for dissipative operators. A coordinate-free approach},  Zap. Nauchn. Sem. Leningrad. Otdel. Mat. Inst. Steklov. (LOMI) 178 (1989), Issled. Linein. Oper. Teorii Funktsii. 18, 57--91, 184--185 (Russian)

English transl.: J. Soviet Math. 61 (1992), no. 2, 1981--2002.

\bibitem[SNF]{SNF} {\sc B. Sz.-Nagy and C. Foia\c s},
{\it Harmonic analysis of operators on Hilbert
space,} Akad\'{e}miai Kiad\'{o}, Budapest, 1970.
 
\end{thebibliography}
